\documentclass[11pt,draft]{article}

\usepackage[T1]{fontenc}
\usepackage{amsfonts,amsmath}

\newtheorem{thm}{theorem}[section]
\newtheorem{theorem}[thm]{Theorem}
\newtheorem{proposition}[thm]{Proposition}
\newtheorem{lemma}[thm]{Lemma}
\newtheorem{corollary}[thm]{Corollary}
\newtheorem{remark}[thm]{Remark}

\newtheorem{definition}[thm]{Definition}

\begin{document}

\title{Graded identities of block-triangular matrices \thanks{partially supported by Fapesp grant No. 2014/10352-4, Fapesp grant No.2014/09310-5, CNPq grant No. 461820/2014-5, CNPq grant No. 480139/2012-1 and CAPES grant No. 99999.001558/2014-05 }}

\author{
Diogo Diniz Pereira da Silva e Silva \thanks{\texttt{diogo@dme.ufcg.edu.br }}
\\
\\
Unidade Acadêmica de Matemática e Estatística\\
Universidade Federal de Campina Grande\\
Campina Grande, PB, Brazil
\\
\\
Thiago Castilho de Mello \thanks{\texttt{tcmello@unifesp.br} }
\\
\\
Instituto de Ci\^encia e Tecnologia\\
Universidade Federal de S\~ao Paulo\\
S\~ao Jos\'e dos Campos, SP, Brazil}

\maketitle

\begin{abstract}

 Let $F$ be an infinite field and $UT(d_1,\dots, d_n)$ be the algebra of
upper block-triangular matrices over $F$. In this paper we describe a basis for the $G$-graded polynomial identities of $UT(d_1,\dots, d_n)$, with an
elementary grading induced by an $n$-tuple of elements of a group $G$ such that the neutral component corresponds to the diagonal of $UT(d_1,\dots,
d_n)$. In particular, we prove that the monomial identities of such algebra follow from the ones of degree up to $2n-1$. Our results generalize for
infinite fields of arbitrary characteristic, previous results in the literature which were obtained for fields of characteristic zero and for
particular $G$-gradings. In the characteristic zero case we also generalize results for the algebra $UT(d_1,\dots, d_n)\otimes C$ with a tensor product grading, where $C$ is a color commutative algebra generating the variety of all color commutative algebras.

\end{abstract}

\section{Introduction}
Lef $F$ be an infinite field and $UT(d_1,\dots, d_n)$ the algebra of upper block triangular matrices. It is the subalgebra of the matrix algebra $M_{d_1+\cdots+d_n}(F)$ consisting of the matrices \[
\left(\begin{array}{cccc}
    A_{11} & A_{12} & \cdots & A_{1n}\\
    0 & A_{22}&\cdots & A_{2n}\\
    \vdots & \vdots & \ddots & \vdots \\
    0&0&\cdots& A_{nn}
\end{array}\right),
\] where $A_{ij}$ is a block of size $d_i\times d_j$. In this paper we  study the graded polynomial identities of upper block triangular matrix algebras $UT(d_1,\dots, d_n)$ over an infinite field $F$. These algebras appear in the classification of minimal varieties (see\cite{GiambrunoZaicev}) and are generalizations of the matrix algebras (when $n=1$) and the algebra $UT_n(F)$ of upper triangular matrices (when $d_1=\dots = d_n=1$).

One of the main problems in the theory of PI-algebras is the (generalized) Specht problem about the existence, for a given class of algebras, of finite basis for the $T$-ideals of identities. This problem for the ordinary identites of associative algebras over a field of characteristic zero was solved by Kemer (see \cite{Kemer}, \cite{Kemer2}). In the case of associative algebras graded by a finite group it was solved by I. Sviridova \cite{Sviridova} in the case of abelian groups and by E. Aljadeff and A. Kanel-Belov \cite{AljadeffBelov} in the general case.
Over fields of positive characteristic however the situations is different and ideals of identities without finite basis exist (see for example \cite{Belov}, \cite{Grishin}, \cite{Shchigolev}). The basis for the graded identities of $UT(d_1,\dots, d_n)$ in our main result (Theorem \ref{Main}) is finite, provided that $G$ is finite.

The algebras of block triangular matrices admit gradings by any group $G$ in which the elementary matrices are homogeneous. These are called \textit{elementary gradings} (or good gradings, see \cite{BarascuDascalescu}). The algebras $UT_n(F)$ admit elementary gradings only (see \cite{ValentiZaicev}). Over an algebraically closed field of characteristic $0$ every grading on $M_n(F)$ by a finite group is obtained by a certain tensor product construction from an elementary grading and a \textit{fine} grading (see \cite{BahturinZaicev}). If moreover the group is abelian an analogous result holds for the algebra $UT(d_1,\dots,d_n)$ (see \cite{ValentiZaicev2}).

Explicit basis for the identities are known for a few algebras only and for the algebras $UT(d_1,\dots, d_n)$ (over an infinite field) the only known basis for the ordinary identities are for the algebras $M_2(F)$ (see \cite{Razmyslov}, \cite{Drensky}, \cite{Koshlukov}) and $UT_n(F)$ (see \cite{Maltsev}). In general, the ideal of identities of $UT(d_1,\dots, d_n)$ is the product of the ideals of identities of the matrix algebras $M_{d_i}$ (see \cite{GiambrunoZaicev0}). An analogous property for the graded identities of block triangular matrix algebras was studied in \cite{DiVincenzoLaScala}. Elementary gradings on $UT_n(F)$ and the corresponding graded identities were studied in \cite{DiVincenzoPlamenValenti} and in particular it was proved that elementary gradings can be distinguished by their graded identities. An analogous result for $UT(d_1,\dots, d_n)$ with an elementary grading by an abelian group was obtained in \cite{DiVincenzoSpinelli}.

When $char F=0$ a complete description of the $\mathbb{Z}_2$-graded identities of $M_2(F)$ (and other PI-algebras) was given in \cite{DiVincenzo}. Analogous basis for the identities of $M_n(F)$ with elementary $\mathbb{Z}$ and $\mathbb{Z}_n$ gradings were determined by Vasilovsky in \cite{Vasilovsky}, \cite{Vasilovsky2}. These results were also established for infinite fields (see \cite{KoshlukovAzevedo}, \cite{Azevedo}, \cite{Azevedo2}) and analogous results were obtained for related algebras (see \cite{DiVincenzo},\cite{DiVincenzoNardozza},\cite{DiVincenzoNardozza2}). Graded identities of $M_n(F)$ were studied more generally in \cite{BahturinDrensky} and in particular a basis for the graded identities of $M_n(F)$ with certain elementary gradings was determined. The elementary $G$-gradings considered are the ones induced by an $n$-tupple of pairwise different elements of $G$. The result considering an infinite field was obtained in \cite{Diogo}. In this case the basis is analogous to the one obtained by Vasilovsky and some monomial identities may be necessary. Recall that a $G$-grading on an algebra $A$ is called \textit{nondegenerate} if the ideal of graded identities of $A$ contains no monomials. These types of gradings were studied in \cite{AljadeffOfir}, \cite{AljadeffOfir2}. Vasilovsky proved that $\mathbb{Z}_n$-grading on $M_n(F)$ is nondegenerate and that for the $\mathbb{Z}$-gradings one needs to consider the monomial identities of degree 1 corresponding to the homogeneous components of dimension 0.

In this paper we prove that the basis given in \cite{BahturinDrensky} holds for the algebras $UT(d_1,\dots,d_n)$ over an infinite field $F$. Moreover
we prove that it is only necessary to include the monomial identities of degree up to $2n-1$ in the basis. In \cite{CentroneMello} a similar result was proved for
$\mathbb{Z}_n$-graded identities. In Section 4 we assume the field $F$ has characteristic zero and generalize the result for the tensor product $UT(d_1,\dots,d_n)\otimes E$, where $E$ denotes the Grassmann algebra with its canonical $\mathbb{Z}_2$-grading, with the tensor product grading. Finally based on the results of \cite{BahturinDrensky}  we generalize our main result to the tensor product $UT(d_1,\dots,d_n)\otimes C$ where $C$ is a color commutative algebra generating the variety of all color commutative algebras. We remark that $C$ is a is a color commutative algebra generating the variety of all color commutative algebras if and only if it has a regular grading. Such gradings were recently studied in \cite{AljadeffOfir}.  

The ideas used in the present paper are similar to those in
 \cite{Azevedo}, \cite{Azevedo2}, \cite{BahturinDrensky}, \cite{CentroneMello}, and \cite{DiVincenzoNardozza3}

\section{Preliminaries}
In this paper we consider associative algebras over an infinite field $F$ and vector spaces are also considered over $F$.

\subsection{Graded algebras and graded polynomial identities}

Let $A$ be an algebra and $G$ a group. A \textit{$G$-grading} on $A$ is a vector space decomposition $A=\oplus_{g\in G}A_g$ compatible with the multiplication of the algebra in the sense that the inclusions \[A_gA_h\subseteq A_{gh}\] hold for any $g$ and $h$ in $G$. A nonzero element $a$ in $\cup_{g\in G}A_g$ is called a \textit{homogeneous element}. Clearly to every homogeneous element $a$ corresponds an element $g$ in $G$ such that $a\in A_g$. We say that this $g$ is the \textit{degree of $a$} in the given $G$-grading. The set $\{g\in G | A_g\neq 0\}$ is the \textit{support} of the grading and is denoted by $supp A$.

A subspace $V$ of $A$ is a \textit{homogeneous subspace} if $V=\oplus_{g\in G}(V\cap A_g)$. A subalgebra $B$ is a \textit{homogeneous subalgebra} if it is homogeneous as a subspace and in this case $B=\oplus_{g\in G}B_g$, where $B_g=B\cap A_g$,  is a $G$-grading on $B$. The $G$-grading on a homogeneous subalgebra $B$ of $A$ is assumed to be this one.

Let $X=\cup X_g$ be a disjoint union of a family of countable sets $X_g=\{x_{g}^{(1)},x_{g}^{(2)},\dots\}$ and $F\langle X|G \rangle$ be the free
associative algebra, freely generated by $X$. If $G$ is clear from the context, we denote it simply by $F\langle X\rangle$. 
A polynomial $f(x_{g_1}^{(1)},\dots, x_{g_n}^{(n)})$ is a graded polynomial identity for $A=\oplus_{g\in G}A_g$ if we have $f(a_{g_1}^{(1)},\dots,
a_{g_n}^{(n)})=0$ for any $a_{g_1}^{(1)}\in A_{g_1}$, \dots, $a_{g_n}^{(n)}\in A_{g_n}$. The set $T_G(A)$ of all graded polynomial identities of $A$
is an ideal of $F\langle X \rangle$ invariant under all graded endomorphisms of this algebra, i.e. it is a $T_G$-ideal. If $S$ is a set of
polynomials in $F\langle X \rangle$ the intersection $U$ of all $T_G$-ideals containing $S$ is a $T_G$-ideal. In this case, we say that $S$ is a
\textit{basis} for $U$. Two sets are equivalent if they generate the same $T_G$-ideal. Since the field $F$ is infinite it is well known that every
polynomial $f$ in $F\langle X \rangle$ is equivalent to a finite collection of multihomogeneous identities. Hence we may reduce our considerations to
multihomogeneous polynomials.

\subsection{Elementary gradings on block-triangular matrices}
Let $(g_1,\dots, g_m)$ be an $m$-tuple of elements of $G$ and $A=M_m(F)$ be the full matrix algebra of order $m$. If we set $A_g$ to be the subspace spanned by the elementary matrices $e_{ij}$ such that $g_i^{-1}g_j=g$ then we have \[A=\oplus_{g\in G}A_g\] and this decomposition is a $G$-grading. Let $B$ be a subalgebra of $A$ generated   by elementary matrices. Then $B$ is a homogeneous subalgebra. In particular $UT(d_1,\dots, d_n)$ is a homogeneous subalgebra of $M_{m}(F)$, where $d_1+\dots+d_n=m$. We say that the $G$-grading on $UT(d_1,\dots, d_m)$ (and more generally on $B$) is the elementary grading induced by $(g_1,\dots, g_n)$.

Let $e$ denote the unit of the group $G$ and consider $B=UT(d_1,\dots, d_m)$ with the elementary grading induced by $(g_1,\dots, g_n)$. The elementary matrices $e_{ii}$ have degree $e$ and therefore the dimension of the component $B_e$ is $\geq n$. We have $dim B_e=n$ if and only if the elements in the $n$-tuple inducing the grading are pairwise distinct. Equivalently the polynomial $x_{e}^{(1)}x_{e}^{(2)}-x_{e}^{(2)}x_{e}^{(1)}$ is a graded identity for $B$.

\subsection{Generic Graded Algebras}

Let $\textbf{g}=(g_1,\dots,g_n)$ be a $n$-tuple of pairwise distinct elements of $G$. Denote by $A=\oplus_{g\in G}A_g$ the algebra $M_n(F)$ with the elementary grading induced by $\textbf{g}$. Let $B$ be a subalgebra of $A$ with basis
$\{e_{i_1j_1},\dots,e_{i_lj_l}\}$ as a vector space. Denote by $G_0$ (resp. $G_0^A$) the support of the grading on $B$ (resp. $A$).

Let $g$ be an element in the support $G_0$ of the grading of $B$. Denote by $D_{\widehat{g}}$ the set of indexes $i \in \{i_1,\dots, i_l\}$ such that for some $j \in \{j_1,\dots,j_l\}$ the matrix unit $e_{ij}$ has degree $g$. Recall that the $n$-tuple $\textbf{g}$ consists of pairwise distinct elements of $G$. This implies that for each $i\in D_{\widehat{g}}$ there exists exactly one index in $\{j_1,\dots, j_l\}$, denoted by $\widehat{g}(i)$, such that $e_{i\widehat{g}(i)}\in B_g$. Thus we obtain a function $\widehat{g}:D_{\widehat{g}}\rightarrow \{j_1,\dots, j_l\}$ for each $g\in G_0$. With this notation $\{e_{i\widehat{g}(i)}|i\in D_{\widehat{g}}\}$ is a basis for $B_g$.

Denote by $\Omega$ the algebra of polynomials in commuting variables \[\Omega=F[\xi_{ij}^{(k)}\,|\,i,j=1,2,\dots,n;\; k=1,2,\dots].\]
The algebra $M_n(\Omega)$ has a natural $G$-grading where the homogeneous component of degree $g$ is the subspace generated by the matrices $m_{ij}e_{ij}$, where $e_{ij}\in A_g$ and $m_{ij}$ is a monomial in $\Omega$.

\begin{definition}\label{def}
For each $g\in G_0$ and each natural number $k$ the element \[\xi^{(k)}_g=\sum_{i\in D_{\widehat{g}}} \xi_{i\widehat{g}(i)}^{(k)}e_{i\widehat{g}(i)}\] of $M_n(\Omega)$  is called a \textbf{graded generic element}. The algebra $G(B)$ generated by the $\xi^{(k)}_g$, $g\in G_0$, $k=1,2,\dots$ is called the \textbf{algebra of graded generic elements of $B$}.
\end{definition}

The algebra $G(B)$ is a homogeneous subalgebra of $M_n(\Omega)$ and is a graded algebra with the inherited grading. If $B=A$ the above construction yields the graded algebra $G(A)$ of generic elements of $A$. The generic element in $G(A)$ corresponding to $g\in G_0^A$ and $k$ will be denoted by $\xi_{g}^{(k,A)}$. The following result is well known.

\begin{theorem}\label{relfree}
Let $F$ be an infinite field. The algebra $G(B)$ is isomorphic as a graded algebra to the relatively free $G$-graded algebra $F\langle X \rangle/ Id_G(B)$.
\end{theorem}
\textit{Proof.}
The homomorphism $\Theta: F\langle X \rangle \rightarrow G(B)$ induced by mapping $x_{g}^{(i)}\mapsto \xi_g^{(i)}$ is clearly onto. Moreover as in the case of the generic matrix algebra (see \cite[Theorem 1.4.4]{GiambrunoZaicev}) we have $ker \Theta = T_G(B)$ and the result follows.
\hfill $\Box$

\section{The main result}

Given $g_1, g_2, \dots, g_p \in G_0$ we consider the composition $\nu = \widehat{g_p}\cdots\widehat{g_1}$ of the corresponding functions. This may not be well defined and we will prove in the next lemma that in this case the monomial
$x_{g_1}^{(1)}\cdots x_{g_p}^{(p)}$ is a graded identity for $B$. Otherwise its domain $D_{\nu}=D_{\widehat{g_p}\cdots\widehat{g_1}}$ is the set of $i\in\{i_1,\dots,i_l\}$ for which the image $\widehat{g_p}(\dots(\widehat{g_1}(i))\dots)$ is well defined. In this case $\{e_{i\nu(i)}|i\in D_{\nu}\}$ is a basis for the subspace spanned by
$B_{g_1}\cdots B_{g_p}$.

\begin{lemma}\label{line}
Let $h_1,h_2,\dots, h_p$  be elements in $G_0$. If $D_{\widehat{h_p}\cdots\widehat{h_1}}=\emptyset$ then $\xi_{h_1}^{(i_1)}\xi_{h_2}^{(i_2)}\cdots \xi_{h_p}^{(i_p)}=0$. Moreover if the set $D_{\widehat{h_p}\cdots\widehat{h_1}}$ is nonempty then the $i$-th line of the matrix  $\xi_{h_1}^{(i_1)}\xi_{h_2}^{(i_2)}\cdots \xi_{h_p}^{(i_p)}$ is
nonzero if and only if
$i\in D_{\widehat{h_p}\cdots\widehat{h_1}}$. In this case if $j=\widehat{h_p}\cdots \widehat{h_1}(i)$, the only nonzero entry in the $i$-th line is a monomial of $\Omega$ in the $j$-th column.
\end{lemma}
\textit{Proof.}
The proof is by induction on the length $p$ of the product. The result for $p=1$ follows directly from Definition \ref{def}. Hence we consider $p>1$ and assume the result for products of length $p-1$. Let us consider first the case
$D_{\widehat{h_{p}}\cdots\widehat{h_1}}\neq\emptyset$. In this case $D_{\widehat{h_{p-1}}\cdots\widehat{h_1}}\neq\emptyset$ and we denote $\nu=\widehat{h_{p-1}}\cdots\widehat{h_1}$. The induction hypothesis implies that there exists monomials $m_i$, where
$i \in D_{\widehat{h_{p-1}}\cdots\widehat{h_1}}$, such that

\begin{equation}\label{a}
\xi_{h_1}^{(i_1)}\xi_{h_2}^{(i_2)}\cdots \xi_{h_p}^{(i_p)}=\left(\sum_{i \in D_{\widehat{h_{p-1}}\cdots\widehat{h_1}}}m_i e_{i \nu(i)}\right) \left(\sum_{j \in D_{\widehat{h_p}}}\xi_{j\widehat{h_p}(j)}^{(i_p)} e_{j\widehat{h_p}(j)}\right).
\end{equation}

Note that $e_{i \nu(i)} e_{j\widehat{h_p}(j)}\neq 0$ for some $j$ if and only if $i \in D_{\widehat{h_{p}}\cdots\widehat{h_1}}$ and in this case the product equals
$e_{i \widehat{h_p}(j)}$. Hence we obtain \[\xi_{h_1}^{(i_1)}\xi_{h_2}^{(i_2)}\cdots \xi_{h_p}^{(i_p)}=\sum_{i\in D_{\widehat{h_{p}}\cdots\widehat{h_1}}} (m_i\xi_{\nu(i)\widehat{h_p}(\nu(i))}^{(i_p)})  e_{i \widehat{h_p}(\nu(i))},\] and the result follows.
Now assume that $D_{\widehat{h_{p}}\cdots\widehat{h_1}}=\emptyset$. If $D_{\widehat{h_{p-1}}\cdots\widehat{h_1}}=\emptyset$ then by the induction hypothesis $\xi_{h_1}^{(i_1)}\xi_{h_2}^{(i_2)}\cdots \xi_{h_{p-1}}^{(i_{p-1})}=0$ and the result holds. Moreover if
$D_{\widehat{h_{p-1}}\cdots\widehat{h_1}}\neq\emptyset$ then we may write the product $\xi_{h_1}^{(i_1)}\xi_{h_2}^{(i_2)}\cdots \xi_{h_p}^{(i_p)}$ as in (\ref{a}). Since $D_{\widehat{h_{p}}\cdots\widehat{h_1}}=\emptyset$ every product $e_{i \nu(i)} e_{j\widehat{h_p}(j)}$ equals zero and therefore  $\xi_{h_1}^{(i_1)}\xi_{h_2}^{(i_2)}\cdots \xi_{h_p}^{(i_p)}=0$.
\hfill $\Box$

\begin{corollary}\label{2}
If a monomial $x_{h_1}^{(i_1)}\dots x_{h_p}^{(i_p)}$ in $F\langle X \rangle$ is a graded identity for $B$ then it is a consequence of a monomial in
$T_G(B)$ of length at most $2n-1$.
\end{corollary}
\textit{Proof.}
The result follows if we prove that every monomial in $T_G(B)$ of length $p>2n-1$ is a consequence of a monomial identity of length at most $p-1$. Let $m=x_{h_1}^{(i_1)}\dots x_{h_p}^{(i_p)}$ be a monomial identity for $B$. Clearly we may assume that $h_i\in G_0$, $i=1,2,\dots,p$. If $D_{\widehat{h_r}\widehat{h_{r-1}}\dots \widehat{h_1}}=\emptyset$ for some $r<p$ then Lema \ref{line} implies that
$\xi_{h_1}^{(i_1)}\dots \xi_{h_r}^{(i_r)}=0$. Hence $x_{h_1}^{(i_1)}\dots x_{h_p}^{(i_p)}$ is an identity for $B$ and $m$ is a consequence of this monomial. Thus we assume that $D_r=D_{\widehat{h_r}\widehat{h_{r-1}}\dots \widehat{h_1}}$ is nonempty for $r<p$ and denote by $I_r$ the image of the composition $\widehat{h_r}\widehat{h_{r-1}}\dots \widehat{h_1}$. Notice that
\begin{equation}\label{incl}
D_1\supseteq D_2 \supseteq \dots \supseteq D_{p-1}\supseteq D_p=\emptyset.
\end{equation}
Assume that there exists $r$ such that $D_r=D_{r+1}=D_{r+2}$. The equality $D_r=D_{r+2}$ implies that $I_r\subseteq
D_{\widehat{h_{r+2}}\widehat{h_{r+1}}}$. Clearly $D_{\widehat{h_{r+2}}\widehat{h_{r+1}}}\subseteq D_{\widehat{h}}$, where $h=h_{r+1}h_{r+2}$.
Therefore $I_r\subseteq D_{\widehat{h}}$ and we conclude that $D_{\widehat{h}\widehat{h_r}\widehat{h_{r-1}}\dots \widehat{h_1}}=D_r$. Since
$D_r=D_{r+2}$ this implies that the compositions $\widehat{h}\widehat{h_r}\widehat{h_{r-1}}\dots \widehat{h_1}$ and
$\widehat{h_{r+2}}\widehat{h_{r+1}}\widehat{h_{r}}\dots \widehat{h_1}$ have the same domain. Moreover the equality in $G$, $h_1h_2\cdots h_{r} h=h_1h_2\cdots
h_{r+1} h_{r+2}$ implies that for every $i\in D_{r+2}$,\linebreak $\widehat{h}\widehat{h_r}\dots \widehat{h_1}(i)=\widehat{h_{r+2}}\widehat{h_{r+1}}\widehat{h_{r}}\dots
\widehat{h_1}(i)$. Hence $\widehat{h}\widehat{h_r}\dots \widehat{h_1}=\widehat{h_{r+2}}\widehat{h_{r+1}}\widehat{h_{r}}\dots
\widehat{h_1}$ and therefore we have $D_{\widehat{h_{p}}\dots \widehat{h_{r+3}}\widehat{h} \widehat{h_r}\dots \widehat{h_1}}=D_p=\emptyset$. It
follows from Lemma \ref{line} that the monomial $m^{\prime}=x_{h_1}^{(i_1)}\dots x_{h_r}^{(i_r)}(x_{h}^{(i_{r+1})})x_{h_{r+3}}^{(i_{r+3})}\dots
x_{h_{p}}^{(i_p)}$, where $i_{r+1}\notin \{i_1,\dots, i_r\}$, is an identity for $B$. Clearly $m$ is a consequence of $m^{\prime}$. It remains only
to verify that if $p>2n-1$ there exists $r$ such that $D_r=D_{r+1}=D_{r+2}$. First notice that if $|D_1|=n$ then
$\{i_1,\dots,i_l\}=\{1,2,\dots,n\}=\{j_1,\dots, j_l\}$ and $\widehat{h_1}$ is a bijection in this set. Therefore $D_p=\emptyset$ implies that
$D_{\widehat{h_p}\widehat{h_{p-1}}\dots \widehat{h_2}}=\emptyset$. By Lemma \ref{line} the monomial $x_{h_2}^{(i_2)}x_{h_3}^{(i_3)}\dots
x_{h_p}^{(i_p)}$ is an identity and clearly $m$ is a consequence of it. Therefore we may assume now that $|D_1|\leq n-1$. In this case there are at
most $n-1$ proper inclusions in (\ref{incl}) and if $p>2n-1$ there are two consecutive equalities, i. e., there exists $r$ such that
$D_r=D_{r+1}=D_{r+2}$. \hfill $\Box$

We consider the following graded polynomials:

\begin{eqnarray}
x_{e}^{(1)}x_{e}^{(2)}-x_{e}^{(2)}x_{e}^{(1)}, \mbox{ if } e\in G_0 \label{(3)}\\
x_{g}^{(1)}x_{g^{-1}}^{(2)}x_{g}^{(3)}-x_{g}^{(3)}x_{g^{-1}}^{(2)}x_{g}^{(1)} \mbox{ if } e \neq g \mbox{ and } B_g\neq 0 \label{(4)} \\
x_{g}^{(1)} \mbox{ if } B_g=0 \label{(5)}.
\end{eqnarray}

\begin{lemma}\label{identities}
The algebra $B$ with the elementary grading induced by an $n$-tuple $g=(g_1,\dots,g_n)$ of pairwise distinct elements of $G$ satisfies the graded polynomial identities $(\ref{(3)}) - (\ref{(5)})$.
\end{lemma}

\textit{Proof.} Clearly the polynomials in $(\ref{(5)})$ are identities for $B$. Since the elements in $\textbf{g}=(g_1,\dots, g_n)$ are pairwise different if $e \in G_0$ the graded generic matrices $\xi_e^{(i)}$ are diagonal. Hence we have the graded identity $(\ref{(3)})$.  Since $(\ref{(4)})$ is multilinear, in order to verify that it is a graded identity substitute $x_{g}^{(1)},x_{g}^{(3)}$ by $e_{i j}, e_{k l} \in B_g$ respectively and $x_{g^{-1}}^{(2)}$ by $e_{r s}\in B_{g^{-1}}$. If $(e_{i j} e_{r s} e_{k l})\neq 0$ then $j=r$ and $s=k$. Moreover $e_{is}$ and $e_{rl}$ are in $A_e$ and therefore $i=s$ and $r=l$. Hence in this case $e_{ij}=e_{kl}$ and the result of the substitution is zero. Analogously if $(e_{kl} e_{r s} e_{ij})\neq 0$ the result is zero. The remaining case to consider is  $(e_{i j} e_{r s} e_{k l})=0=(e_{kl} e_{r s} e_{ij})$ and the result is also $0$. \hfill $\Box$

\begin{proposition}\cite[Lemma 4.6]{Diogo}\label{J}
Let $U_A$ denote the $T_G$-ideal generated by the identities $(\ref{(3)})-(\ref{(5)})$ satisfied by the matrix algebra $A$ and let $\xi_g^{(i,A)}$, $g\in G_0^A$, $i=1,2,\dots$ denote the generic elements in $G(A)$. If the monomials $m(x_{h_1}^{(1)},\dots,x_{h_p}^{(p)})$ and $n(x_{h_1}^{(1)},\dots,x_{h_p}^{(p)})$ in $F\langle X \rangle$ are such that the matrices $n(\xi_{h_1}^{(1,A)}\dots,\xi_{h_p}^{(p,A)})$ and $n(\xi_{h_1}^{(1,A)}\dots,\xi_{h_p}^{(p,A)})$ have the same position the same non-zero entry then \[m(x_{h_1}^{(1)},\dots,x_{h_p}^{(p)})\equiv n(x_{h_1}^{(1)},\dots,x_{h_p}^{(p)}) \mbox{ modulo }U_A.\]
\end{proposition}

Next we generalize this proposition to the case of a subalgebra $B$ of $A=M_n(F)$ generated by elementary matrices. Note that the algebra $G(B)$ is a homomorphic image of the algebra $G(A)$ by Theorem \ref{relfree}. The homomorphism constructed in the following remark will be usefull.

\begin{remark}\label{homomorphism}
We construct a homomorphism from $G(A)$ to $G(B)$ as follows: the map $x_{ij}^{(k)}\mapsto \chi_{ij}x_{ij}^k$ where $\chi_{ij}=1$ if $e_{ij}\in B_g$ and $\chi_{ij}=0$ if $e_{ij}\notin B_g$ induces an endomorphism $\theta$ of $\Omega$ extending this map. Hence $\Theta:M_n(\Omega)\rightarrow M_n(\Omega)$ given by $\Theta(\sum p_{ij}e_{ij})=\sum \theta(p_{ij})e_{ij}$ is an endomorphism of $M_n(\Omega)$. From the definiton of $\theta$ if follows that
$\Theta (\xi_{ij}^{(k,A)})=\xi_{ij}^{(k)}$ and therefore $\Theta (G(A))=G(B)$. The restriction to $G(A)$ gives the desired homomorphism (also denoted by $\Theta$).
\end{remark}

\begin{corollary}\label{c}
Let $B$ be a subalgebra of $M_n(F)$ generated   by elementary matrices with the induced $G$-grading and $U_0$ be the $T_G$-ideal generated
by the identities $(\ref{(3)})-(\ref{(5)})$ satisfied by the graded algebra $B$. If $m(x_{h_1}^{(1)},\dots,x_{h_p}^{(p)})$ and $n(x_{h_1}^{(1)},\dots,x_{h_p}^{(p)})$ are two
monomials in $F\langle X \rangle$ such that the matrices $m(\xi_{h_1}^{(1)},\cdots,\xi_{h_p}^{(p)})$ and $n(\xi_{h_1}^{(1)},\cdots,\xi_{h_p}^{(p)})$ have in
the same position the same nonzero entry then
\[m(x_{h_1}^{(1)},\dots,x_{h_p}^{(p)})\equiv n(x_{h_1}^{(1)},\dots,x_{h_p}^{(p)}) \mbox{ modulo }U_0.\]
\end{corollary}

\textit{Proof.} Let $\tilde{m} (x_{g_1}^{(1)}\cdots x_{g_n}^{(n)})$ be a monomial in $F\langle X \rangle$. Let $\Theta$ be the homomorphism constructed in the previous remark. We have
\begin{equation}\label{eq}
 \Theta (\tilde{m}(\xi_{g_1}^{(i_1,A)}\cdots \xi_{g_1}^{(i_n,A)})) =  \tilde{m}(\xi_{g_1}^{(i_1)}\cdots \xi_{g_n}^{(i_n)}).
\end{equation}
It follows from Lemma \ref{line} that the entries of $\tilde{m}(\xi_{g_1}^{(i_1,A)}\cdots \xi_{g_1}^{(i_1,A)})$ are monomials in
$\Omega$. Note that if $p$ is a monomial in $\Omega$ then $\theta (p)$ is either $0$ or $p$. Hence (\ref{eq}) implies that the nonzero entries of $\tilde{m}(\xi_{g_1}^{(i_1)}\cdots \xi_{g_n}^{(i_n)})$ equal the corresponding entries of $\tilde{m}(\xi_{g_1}^{(i_1,A)}\cdots \xi_{g_n}^{(i_n,A)})$. Thus the matrices $m(\xi_{h_1}^{(1,A)},\cdots,\xi_{h_p}^{(p,A)})$ and $n(\xi_{h_1}^{(1,A)},\cdots,\xi_{h_p}^{(p,A)})$ have in the same position the same nonzero entry. Therefore Proposition \ref{J} imples that $m\equiv n$ modulo $U_A$. The result is then a consequence of the inclusion $U_A\subseteq U_0$. To prove this we verify that every generator of $U_A$ is in $U_0$ and this follows from the inclusion $G_0\subset G_0^A$.
\hfill $\Box$

\begin{theorem}\label{Main}
Let $G$ be a group and let $\textbf{g}=(g_1,\dots, g_n)\in G^{n}$ induce an elementary $G$-grading of $M_{n}(F)$, where the elements $g_1,\dots, g_n$ are pairwise different. If $B$ is a subalgebra of $M_n(F)$ generated   by elementary matrices $e_{ij}$ then a basis of the graded polynomial identities of $B$ consists of $(\ref{(3)})-(\ref{(5)})$ and a finite number of identities of the form
$x_{h_1 1}\dots x_{h_p p}$, where $2\leq p\leq 2n-1$.
\end{theorem}
\textit{Proof.} Let $U$ be the $T_G$-ideal of $F\langle X \rangle$ generated by the polynomials (\ref{(3)}) -- (\ref{(5)})
together with the monomial identities $x_{h_1}^{(1)}\dots x_{h_p}^{(p)}$ of $B$ with $2\leq p\leq 2n-1$. It follows from Lemma \ref{identities} that $U\subseteq T_G(B)$. Hence to prove the theorem it is enough to show that every multihomogeneous $G$-graded identity of $B$ lies in $U$. Assume, on the contrary, that $f$ is a multihomogeneous graded identity that does not lie in $U$. We write $f\equiv\sum_{i=1}^k \alpha_im_i$  modulo $U$, where the $\alpha_i$ are non-zero scalars and $m_i$ are monomials in $F\langle X \rangle$. We may assume that the number $k$ of nonzero coefficients
is minimal. If $k=1$ then $m_1$ is an identity for $B$ and Corollary \ref{2} implies that it lies in $U$ which is a contradiction. We now consider
$k>1$. Denote by $\overline{m_i}$ the matrix in $G(B)$ that is the result of substituting every variable $x_{g}^{(j)}$ in $m_i$ for the corresponding
generic matrix $\xi_{g}^{(i)}$. By the minimality of $k$ the monomials $m_i$ are not identities for $B$ and in particular $\overline{m_1}$ has a nonzero
entry. Moreover we have
\[-\alpha_1\overline{m_1}=\sum_{i=2}^k\alpha_i\overline{m_i}.\] It follows from Lema \ref{line} that the nonzero entries of the matrices $\overline{m_i}$ are monomials in $\Omega$. Therefore there exists a $j>1$ such that $\overline{m_j}$ and $\overline{m_1}$ have in the same position the same nonzero entry. Thus Corollary \ref{c} implies that $m_1\equiv m_j$ modulo $U$.  Hence $f\equiv(\alpha_1+\alpha_j)m_1+\sum_{i\neq j}m_i$ modulo $U$.
This last polynomial is an identity for $B$ that does not lie in $U$ with fewer nonzero coefficients than $f$ and this is a contradiction.
\hfill $\Box$

Recall that a $G$-grading on an algebra $A$ is called nondegenerate if for every integer $r$ and any tuple $(g_1,\dots, g_r)\in G^r$ the monomial $x_{g_1}^{(1)}\cdots x_{g_r}^{(r)}$ is not a graded identity for $A$ (see \cite[Observation 2.2]{AljadeffOfir}). A stroger condition is that $A_gA_h=A_{gh}$ for every $g,h\in G$ and in this case the grading is called \textit{strong}. The $\mathbb{Z}_n$-grading considered by Vasilovsky in \cite{Vasilovsky2} is strong and in particular nondegenerate and the basis determined consists of $(\ref{(3)})$ and $(\ref{(4)})$. In the next corollary we consider elementary $G$-gradings on $M_n(F)$ that are closely related to this grading.

\begin{corollary}
Let $G$ be a finite group with unit $e$ and let $M_n(F)$ be endowed with an elementary grading such that $x_e^{(1)}x_e^{(2)}-x_e^{(2)}x_e^{(1)}$ is a graded polynomial identity. If the $G$-grading is nondegenerate then a basis for the graded identities of $M_n(F)$ consists of the polynomials $(\ref{(3)})$
and $(\ref{(4)})$. Moreover in this case the grading is strong and $G$ is a group of order $n$.
\end{corollary}

\textit{Proof.} Let $\textbf{g}=(g_1,\dots,g_n)\in G^n$ be a tuple inducing the elementary grading. If $g_i= g_j$ for some $i\neq j$ then the
elementary matrices $e_{ij}$ and $e_{ji}$ have degree $e$ and $e_{ij}e_{ji}-e_{ji}e_{ij}\neq 0$. Hence $x_e^{(1)}x_e^{(2)}-x_e^{(2)}x_e^{(1)}$ is not a graded
identity and this is a contradiction. Thus the elements in the tuple $\textbf{g}$ are pairwise different. Since the grading is nondegenerate it
follows from Theorem \ref{Main} that the polynomials $(\ref{(3)})$ and $(\ref{(4)})$ are a basis for the graded identities. Now we prove the last assertion. Note
that since $\textbf{g}$ consists of pairwise different elements we have $|G|\geq n$. We claim that if $|G|>n$ then there exists $g_1,\dots, g_n \in
G$ such that $x_{g_1}^{(1)}\cdots x_{g_n}^{(n)}$ is a graded identity. Clearly it follows from this claim that $|G|=n$. We construct the sequence as follows:
since $|G|> n$ we let $g_1\in G$ such that none of the elementary matrices $e_{11},e_{12},\dots,e_{1n}$ have degree $g_1$. Clearly the first line of
$\xi_{g_1}^{(1)}$ is zero. Then we choose $g_2$ such that the second line of $\xi_{g_1g_2}^{(1)}$ is zero. Inductively we choose $g_i$ such that the $i$-th
line of $\xi_{g}^1$ is zero, where $g=g_1\cdots g_{i}$. Note that the first line of $\xi_{g_1}^{(1)}\xi_{g_2}^{(2)}$ is zero since the first line of
$\xi_{g_1}^{(1)}$ is zero. Moreover the second line of $\xi_{g_1}^{(1)}\xi_{g_2}^{(2)}$ is also zero because the second line of $\xi_{g_1g_2}^{(1)}$ is zero. It
follows by induction that the first $i$ lines of $\xi_{g_1}^{(1)}\cdots \xi_{g_i}^{(i)}$ are zero. Hence $\xi_{g_1}^{(1)}\cdots \xi_{g_n}^{(n)}=0$ and it follows
from Lemma \ref{line} that $x_{g_1}^{(1)}\cdots x_{g_n}^{(n)}$ is a graded identity for $M_n(F)$. Now we prove that the grading is strong. Let $g\in G$. Note
that for each $i$ the elementary matrices $e_{i1},\dots, e_{in}$ have pairwise different degrees and since $|G|=n$ the sequence of degrees is just a
reordering of the elements of $G$. Thus there exists $j$ such that $e_{ij}$ has degree $g$. Hence we obtain $A_gA_h=A_{gh}$ for any $g,h\in G$.
\hfill $\Box$

\begin{remark}
The proof of the last assertion in the previous lemma is based on the proof of Lemma 3.3 in \cite{AljadeffOfir}.  In this lemma a characterization of
nondegenerate gradings on finite dimensional $G$-simple algebras is given.
\end{remark}

\begin{corollary}
Let $G$ be a group. If $UT(d_1,\dots, d_n)$ has an elementary grading such that the polynomials $(\ref{(3)})$ and $(\ref{(4)})$ are a basis for the graded polynomial
identities of this graded algebra then $n=1$,i.e., $UT(d_1,\dots, d_n)=M_{d_1}(F)$.
\end{corollary}
\textit{Proof.}
If $n>1$ then we apply the previous lemma to each block $A_{ii}$ to obtain a monomial that is a graded identity for $M_{d_i}$. The product of copies of these monomials in disjoint sets of variables is a monomial $m$ such that the result of any substitution lies in the jacobson radical $J$ of $UT(d_1,\dots, d_n)$. Since $J$ is a nilpotent ideal, say $J^k=0$, the product of $k$ copies of $m$ in disjoint sets of variables is a monomial identity.
\hfill $\Box$

\section{Matrices over the Grassmann algebra}

We now turn our attention to matrices over the Grassmann algebra.

In this section we suppose $F$ is a field of characteristic zero and we denote by $E$ the Grassmann algebra of an infinite dimensional vector space
over $F$ with its natural $\mathbb{Z}_2$-grading $E=E_0\oplus E_1$ induced by the length of its monomials. For more information concerning the
Grassmann algebra, see \cite{DrenskyBook} .

We use the results of the previous sections and results of \cite{DiVincenzoNardozza3} to find a basis for the $G\times \mathbb{Z}_2$-graded
polynomial identities of $UT(d_1,\dots,d_n; E)$: the algebra of block-triangular matrices over the Grassmann algebra, which is isomorphic to the
tensor product $UT(d_1,\dots,d_n) \otimes E$, and more generally of the algebra $B\otimes E$, where $B$ is a $G$-graded subalgebra of $M_n(F)$
generated by elementary matrices with an elementary grading induced by an $n$-tuple $(g_1,\dots,g_n)$ of pairwise distinct elements of $G$.

If $B$ is a $G$-graded algebra the algebra $B\otimes E$ has a natural $G\times \mathbb{Z}_2$-grading induced by the gradings of $B$ and of $E$. In
such grading, the homogeneous component of degree $(g,\delta)$ is $(B\otimes E)_{(g,\delta)}=B_g\otimes E_\delta$.

In order to work with the $G\times \mathbb{Z}_2$-graded identities of $B\otimes E$, we now consider the free associative algebra $F\langle Z
\rangle$, with $Z=X'\cup Y'$, where $X'=\cup X'_g$ is the set of graded variables with $G\times \mathbb{Z}_2$-degree $(g,0)$ and $Y'=\cup Y'_g$ is
the set of graded variables with $G\times \mathbb{Z}_2$-degree $(g,1)$. We denote elements of $X'_g$ and $Y'_g$ respectively by $x_{g}^{(i)}$ and
$y_{g}^{(i)}$, for $i\in \mathbb{N}$ and $g\in G$. From now on, the variables labeled as $z_{g}^{(i)}$ may be $x_{g}^{(i)}$ or $y_{g}^{(i)}$.

Recall that in \cite{DiVincenzoNardozza3} the authors define a map $\zeta_J$, for $J\subseteq \mathbb{N}$, which maps multilinear identities of $B$
into identities of $B\otimes E$. Such map is defined as follows.

First, we observe that $F\langle Z\rangle $ is both a $\mathbb{Z}_2$-graded algebra and $G$-graded algebra.  Concerning the $\mathbb{Z}_2$-grading of $F\langle Z\rangle $, one defines the map $\zeta$ as follows. If $m$ is a multilinear monomial let
$i_1<\cdots < i_k$ be the indices with odd $\mathbb{Z}_2$-degree occurring in $m$. Then for some $\sigma\in \rm{Sym}(\{i_1, \dots,i_k\})$, we write
\[m=m_0y_{g_1\sigma(i_1)}m_1\cdots y_{g_k\sigma(i_k)}m_{k+1}\]
where $m_0,\dots,m_{k+1}$ are monomials on even variables only. Then we define
\[\zeta(m)=(-1)^\sigma m\]

\begin{definition}
Let $J\subseteq \mathbb{N}$. We define $\varphi_J: F\langle X\rangle \longrightarrow F\langle Z\rangle$ to be the unique $G$-homomorphism of algebras
defined by

\[\varphi_J(x_{g}^{(i)})=\left\{\begin{array}{cc}
   x_{g}^{(i)} & \text{ if } i\not \in J \\
   y_{g}^{(i)} & \text{ if } i \in J
 \end{array}\right.\]
Also for a multilinear monomial $m$ we define $\zeta_J(m)=\zeta(\varphi_J(m))$
\end{definition}

The map $\zeta_J$ extends by linearity to the space of all multilinear polynomials in $F\langle X\rangle$ and for each multilinear polynomial in
$F\langle X\rangle$,  $\zeta_J(f)$ is also a multilinear polynomial in $F\langle Z\rangle$.

We now recall Theorem 11 of \cite{DiVincenzoNardozza3}.
\begin{theorem}\label{Theorem 11}
Let $A$ be a $G$-graded algebra and $\mathcal{E}\subset F\langle X|G \rangle$ be a system of multilinear generators for $T_G(A)$. Then the set
 \[ \{\zeta_J(f)\,|\, f\in \mathcal{E}, \; J\subseteq \mathbb{N}\}\] is system of multilinear generators of $T_{G\times \mathbb{Z}_2}(A\otimes E)$
\end{theorem}

Since the basis of the graded polynomial identities of $UT(d_1,\dots,d_m)$, described in Theorem \ref{Main} contains polynomials in at most $2n-1$
variables, it is enough to consider $J\subset \{1,\dots,2n-1\}$.

\begin{lemma} Applying the map $\zeta_J$ to the polynomial $x_{e}^{(1)}x_{e}^{(2)}-x_{e}^{(2)}x_{e}^{(1)}$ we obtain up to endomorphisms of $F\langle X|G\times
\mathbb{Z}_2\rangle$ the following polynomials

\begin{eqnarray}
              x_{e}^{(1)}z_{e}^{(2)}-z_{e}^{(2)}x_{e}^{(1)}\label{7}\\
              y_{e}^{(1)}y_{e}^{(2)}+y_{e}^{(2)}y_{e}^{(1)}
\end{eqnarray}
Applying the map $\zeta_J$ (for $J\subseteq \{1,2,3\}$) to the polynomial
$x_{g}^{(1)}x_{g^{-1}}^{(2)}x_{g}^{(3)}-x_{g}^{(3)}x_{g^{-1}}^{(2)}x_{g}^{(1)}$ we obtain up to endomorphisms of $F\langle X|G\times
\mathbb{Z}_2\rangle$ the polynomials
\begin{eqnarray}
              z_{g}^{(1)}x_{g^{-1}}^{(2)}x_{g}^{(3)}-x_{g}^{(3)}x_{g^{-1}}^{(2)}z_{g}^{(1)} \\
              x_{g}^{(1)}y_{g^{-1}}^{(2)}x_{g}^{(3)}-x_{g}^{(3)}y_{g^{-1}}^{(2)}x_{g}^{(1)} \\
              y_{g}^{(1)}y_{g^{-1}}^{(2)}x_{g}^{(3)}+x_{g}^{(3)}y_{g^{-1}}^{(2)}y_{g}^{(1)} \\
              y_{g}^{(1)}x_{g^{-1}}^{(2)}y_{g}^{(3)}+y_{g}^{(3)}x_{g^{-1}}^{(2)}y_{g}^{(1)} \\
              y_{g}^{(1)}y_{g^{-1}}^{(2)}y_{g}^{(3)}+y_{g}^{(3)}y_{g^{-1}}^{(2)}y_{g}^{(1)}\label{13}
\end{eqnarray}
Finally, if $m=x_{g_1}^{(1)}\cdots x_{g_p}^{(p)}$ is a $G$-graded monomial identity of the algebra $B$, generated by elementary matrices, for some
$1\leq p\leq 2n-1$ then up to some endomorphism of $F\langle X|G\times \mathbb{Z}_2\rangle$, $\zeta_J(m)=z_{g}^{(1)}\cdots z_{g}^{(p)}$
\end{lemma}
\textit{Proof.} The proof consist of several applications of the map $\zeta_J$ for $J\subseteq \mathbb{N}$. For
$f_1=x_{e}^{(1)}x_{e}^{(2)}-x_{e}^{(2)}x_{e}^{(1)}$, it is enough to consider $J\subseteq \{1,2\}$. So consider $J=\{1\}$ and $J=\{2\}$. Then
$\zeta_{\{1\}}(f_1)=y_{e}^{(1)}x_{e}^{(2)}-x_{e}^{(2)}y_{e}^{(1)}$ and $\zeta_{\{2\}}(f_1)=-(y_{e}^{(2)}x_{e}^{(1)}-x_{e}^{(1)}y_{e}^{(2)})$, and the
latter is the image of the former, by the endomorphism of $F\langle Z\rangle$, which permutes the indexes 1 and 2 of the variables and multiply the
result by $-1$. For this reason, up to an endomorphism of $F\langle Z\rangle$, the image of $f_1$ is
$x_{e}^{(1)}z_{e}^{(2)}-z_{e}^{(2)}x_{e}^{(1)}$, for $|J|\leq 1$. If $J=\{1,2\}$, one obtains $\zeta_{J}=y_{e}^{(1)}y_{e}^{(2)}+y_{e}^{(2)}y_{e}^{(1)}$.

Similarly, one obtains the images  by $\zeta_J$ of the polynomials \[f_2=
x_{g}^{(1)}x_{g^{-1}}^{(2)}x_{g}^{(3)}-x_{g}^{(3)}x_{g^{-1}}^{(2)}x_{g}^{(1)} \text{ and } m=x_{g_1}^{(1)}\cdots x_{g_p}^{(p)}\] \hfill $\Box$

By applying the above lemma and theorem we obtain:

\begin{corollary}
Let $F$ be a field of characteristic zero, $G$ be a group and let $\textbf{g}=(g_1,\dots, g_n)\in G^{n}$ induce an elementary $G$-grading of
$M_{n}(F)$, where the elements $g_1,\dots, g_n$ are pairwise different. If $B$ is a subalgebra of $M_n(F)$ generated   by matrix units $e_{ij}$, then
a basis of the graded polynomial identities of the algebra $B\otimes E$ consists of the polynomials $(\ref{7})-(\ref{13})$ and a finite number of identities of
the form $z_{g_1}^{(1)}\dots z_{g_p}^{(p)}$, with $2\leq p\leq 2n-1$, for each $g_1$,\dots,$g_p\in G_0$ such that $x_{g_1}^{(1)}\dots x_{g_p}^{(p)}$
is a graded identity of $B$.
\end{corollary}

\begin{remark}
It is interesting to observe that in characteristic $p$ case the map $\zeta_J$ also maps multilinear identities of $B$ into multilinear identities of
$B\otimes E$. But such identities may not be enough to generate all $G\times \mathbb{Z}_2$-graded identities of $B\otimes E$, since in positive
characteristic, the identities may not be generated by the multilinear ones.

As an example one can consider the field $F$, as the algebraic closure of the prime field $\mathbb{Z}_p$, graded by the trivial group. The ideal of
identities of $F$ are generated by the polynomial $[x_1,x_2]$, but the algebra $F\otimes E$, which is isomorphic to $E$, satisfies the
$\mathbb{Z}_2$-graded identity
\[St_p(y_1,\dots,y_p)=\sum_{\sigma\in S_p}(-1)^\sigma y_{\sigma(1)}\cdots y_{\sigma(p)}\]
which is not in the $T_{\mathbb{Z}_2}$-ideal generated by the image of $[x_1,x_2]$ by $\zeta_J$.

Problems involving relations between identities in positive characteristic and in characteristic zero are quite difficult. See for example
\cite[Problem e), p. 185]{Procesi}.
\end{remark}

\section{Color Commutative Superalgebras}

In this section we study the connection between identities of a $G$-graded algebra $R$ and the identities of the tensor product of $R$ by an $H$-graded color commutative superalgebra, $C$, where $H$ is an abelian group. Again, we work over a field of characteristic zero.

If $H$ is an abelian group, written additively, let $\beta:H\times H \longrightarrow F^*$ be a skew-symmetric bicharacter, i.e., a function satisfying for all $g, h, k\in H$, the following properties

\[\begin{array}{c}
\beta(g+h,k)=\beta(g,k)\beta(h,k) \\
\beta(g,h+k)=\beta(g,h)\beta(g,k) \\
\beta(g,h)=\beta(h,g)^{-1}
\end{array}\]

Now, if $C=\oplus_{h\in H}C_h$, we define the $\beta$-commutator 

\[ [a,b]_\beta=ab-\beta(h,k)ba \]
for $a\in C_h$ and $b\in C_k$ and extend it by linearity to $C$. We say that $C$ is $\beta$-commutative if $[a,b]_\beta=0$, for all $a$, $b \in C$. If $\beta$ is fixed we call $\beta$-commutative algebras ``color commutative superalgebras'' \cite{BaMiPeZa}.

As examples, if one considers $\beta\equiv 1$, $\beta$-commutative algebras are simply commutative algebras. If one considers the Grassmann algebra $E$, with its usual $\mathbb{Z}_2$-grading, one can see that defining $\beta(g,h)=1$, if $g=0$ or $h=0$ and $\beta(g,h)=-1$ otherwise, one obtains that $[x,y]_\beta=0$, for all $x$, $y\in E$, i.e., $E$ is a color commutative superalgebra.

The main result of this section is a generalization of Theorem \ref{Theorem 11} above \cite[Theorem 11]{DiVincenzoNardozza3}, i.e., we want to replace $E$ by an arbitrary $H$-graded color commutative superalgebra $C$ in the above theorem.

If $R$ is a $G$-graded algebra and $C$ is a $H$-graded color commutative superalgebra, of course, $R\otimes C$ is a $G\times H$-graded algebra.

If $G$ is a group and ${\bf g}=(g_1,\dots,g_n)\in G^n$, we denote by $P_{\bf g}$, the subspace of $F\langle X|G\rangle$ generated by
$\{x_{g_{\sigma(1)}\sigma(1)}\cdots x_{g_{\sigma(n)}\sigma(n)}\,|\, \sigma \in S_n\}$.
The space of multilinear polynomials of degree $n$ of $F\langle X\,|\, G\rangle$ is defined by $P_n^G=\bigoplus_{{\bf g}\in G^n} P_{\bf g}$. It is
well known that if $F$ has characteristic zero, the $G$-graded identities of a $G$-graded $F$-algebra $A$, follow from the ones in $P_{\bf g}$, for ${\bf g}\in
G^n$ and $n\in \mathbb{N}$.

If ${\bf g}=(g_1,\dots,g_n)\in G^n$ and ${\bf h}=(h_1,\dots,h_n) \in H^n$, we denote ${\bf g\times h}=((g_1,h_1),\dots,(g_n,h_n))\in (G\times H)^n$.

For each sequence of elements of $H$, ${\bf h}=(h_i)_{i\in \mathbb{N}}$, we define a map $\varphi_{\bf h}:F\langle X|G \rangle\longrightarrow F\langle X|G\times H \rangle$ as the unique homomorphism of $G$-graded algebras satisfying $\varphi_{\bf h}(x_{g_ii})=x_{(g_i,h_i)i}$.  Here the $G$-grading on $F\langle X| G\times H\rangle$ is the one induced by $\deg_G(x_{(g,h)i})=g$. If ${\bf g}=(g_1,\dots,g_n)\in G^n$ we denote  by the same symbol, ${\bf g}$, the sequence $(g_i)_{i\in \mathbb{N}}$ such that $g_i$ is the neutral element of $G$, if $i> n$.

Now, for each multilinear monomial $m\in F\langle X|G\times H \rangle$, in the variables $x_{(g_{i_1},h_{i_1})i_1}$, \dots, $x_{(g_{i_k},h_{i_k})i_k}$, with $i_1<\cdots<i_k$, we may write $m=x_{(g_{i_{\sigma(1)}}h_{i_{\sigma(1)}})i_{\sigma(1)}}\cdots x_{(g_{i_{\sigma(k)}}h_{i_{\sigma(k)}})i_{\sigma(k)}}$, for some permutation $\sigma$.

If in the free $H$-graded color commutative superalgebra, we have \[x_{h_{i_1}i_1}\cdots x_{h_{i_k}i_k}=\lambda_{{\bf h}\sigma}x_{h_{i_{\sigma(1)}}i_{\sigma(1)}}\cdots x_{h_{i_{\sigma(k)}}i_{\sigma(k)}},\] with $\lambda_{{\bf h} \sigma}\in F^*$, we define $\zeta(m)=\lambda_{{\bf h}\sigma}m$.

Also, for each multilinear monomial $m\in P_n^G$, we define the map $\phi_{\bf h}$ on $m$ as  $\phi_{\bf h}(m)=\zeta(\varphi_h(m))$ and extend $\phi_{\bf h}$ to $P_n^G$ by linearity.

\begin{theorem}\label{3.1}\cite[Theorem 3.1]{BahturinDrensky}
Let $C$ be an $H$-graded color commutative superalgebra, generating the variety of all $H$-graded color commutative superalgebras and let $R$ be any
$G$-graded algebra. If $f(x_{g_11},\dots,x_{g_nn})$ is a multilinear $G$-graded polynomial and ${\bf h}=(h_1,\dots,h_n)\in H^n$, then
$f(x_{g_11},\dots,x_{g_nn})=0$ is a graded polynomial identity for the $G$-graded algebra $R$ if and only if
$\phi_{\bf{h}}(f)(x_{(g_1,h_1)1},\dots,x_{(g_n,h_n)n})=0$ is a graded polynomial identity for the $(G\times H)$-graded algebra $R\otimes C$.
\end{theorem}

Although the above result associates $G$-graded identities  of $R$ with $G\times H$-identities of $R\otimes C$, it is not enough for our purposes. We want to construct a basis for the $G\times H$-graded identities of $R\otimes C$ starting from a basis of the $G$-graded identities of $R$. In order to construct such basis, we need to generalize some lemmas used in the proof of Theorem \ref{Theorem 11} \cite[Theorem 11]{DiVincenzoNardozza3}, namely Lemmas 5, 9 and 10 of \cite{DiVincenzoNardozza3} to the case of tensor product by color commutative superalgebras. This follows below.

We remark that the set $J$ and the map $\zeta_J$ in section 4 plays the same role of the sequence {$\bf h$} and the map $\phi_{\bf h}$ in this section.

\begin{lemma}\label{5}
Let ${\bf g}\in G^n$ and ${\bf h}\in H^n$. If $f\in T_{G\times H}(R\otimes C)\cap P_{\bf g\times h}$. Then there exists $f_0\in P_{\bf g}\cap T_G(R)$
such that $f=\phi_{\bf h}(f_0)$.
\end{lemma}
\textit{Proof.} If $f\in T_{G\times H}(R\otimes C)\cap P_{\bf g\times h}$ we may write \[f=\sum_{\sigma\in S_n} \alpha_\sigma
x_{(g_{\sigma(1)},h_{\sigma(1)})\sigma(1)}\cdots x_{(g_{\sigma(n)},h_{\sigma(n)})\sigma(n)},\] for some $\alpha_\sigma\in F$. Now we define $f_0$ as
\[f_0=\sum_{\sigma\in S_n} \alpha_\sigma\lambda_{{\bf h}\sigma}^{-1}
x_{g_{\sigma(1)}\sigma(1)}\cdots x_{g_{\sigma(n)}\sigma(n)},\] where $\lambda_{{\bf h}\sigma}\in F^*$ is the coefficient used in the definition of $\phi_{\bf h}$.

Now it is easy to observe that $\phi_{\bf h}(f_0)=f$. By Theorem \ref{3.1}, $f_0 \in T_G(R)$ if and only if $f\in T_{G\times H}(R\otimes C)$ and this
proves the lemma. \hfill $\Box$

\begin{lemma}\label{9}
Let $u_1$, \dots, $u_m$ be monomials in $F\langle X|G\rangle$ such that for some $n\geq m$, ${\bf g}\in G^n$ and ${\bf h}\in H^n$, $u=u_1\cdots u_m \in P_{\bf g\times h}$. Consider
${\bf h'}=(h_1',\dots, h_m')\in H^m$ such that $h_i'=\deg_H(u_i)$. Then,  there exists $\gamma\in F^*$ such that
for every $\sigma\in S_m$,
\[\zeta(u_{\sigma(1)}\cdots u_{\sigma(m)})=\gamma\lambda_{{\bf h'}\sigma}u_{\sigma(1)}\cdots u_{\sigma(m)}\]
\end{lemma}
\textit{Proof.}  By the definition of the map $\zeta$, there exist $\gamma\in F^*$ such that
\[\zeta(u_1\cdots u_m)=\gamma u_1\cdots u_m.\]
Since in the free $H$-graded color commutative superalgebra we have
\[\lambda_{{\bf h'}\sigma}x_{h_{\sigma(1)}'\sigma(1)}\cdots x_{h_{\sigma(m)}'\sigma(m)}=x_{h_1'1}\cdots x_{h_m'm}\]
and for each $i\in \{1,\dots,m\}$,  $\deg(u_{i})=h_i'=\deg(x_{h_{i}'i})$, we obtain that
\[\zeta(u_{\sigma(1)}\cdots u_{\sigma(m)})=\gamma \lambda_{{\bf h'}\sigma}u_{\sigma(1)}\cdots u_{\sigma(m)},\] which proves the lemma. \hfill $\Box$

\begin{lemma}\label{10}
Let $f=f(x_{g_11},\dots,x_{g_mm})\in P_m^G$, and  $w_1,\dots,w_m$ monomials in $F\langle X|G\rangle$ such that $\deg_G(w_i)=g_i$, for each $i$ and $f(w_1,\dots,w_m)\in P_n^{G}$. If ${\bf h}\in H^n$, there exist ${\bf h'}\in H^m$ and homogeneous elements $b_1,\dots, b_m\in F\langle X|G\times H\rangle$ such that
$\phi_{\bf h}(f(w_1,\dots,w_m))=\gamma
\phi_{\bf h'}(f)(b_1,\dots,b_m)$
\end{lemma}

\textit{Proof.} For each $i\in\{1,\dots,m\}$, define $b_i:=\varphi_{\bf h}(w_i)$ and $h_i':=\deg_H(b_i)$ and let ${\bf h'}:=(h_1',\dots,h_m')\in H^m$.

Suppose now that \[f(x_{g_11},\dots,x_{g_mm})=\sum_{\sigma\in S_m}c_\sigma x_{g_{\sigma(1)}\sigma(1)}\cdots x_{g_{\sigma(m)}\sigma(m)}.\] Then, if $\tilde{f}=\phi_{\bf h'}(f)$, we have

\[ \tilde{f}(x_{(g_1,h_1')1},\dots,x_{(g_m,h_m')m})=\hspace{-.03cm}\sum_{\sigma\in S_m}c_\sigma \lambda_{{\bf h'}\sigma} x_{(g_{\sigma(1)},h_{\sigma(1)}')\sigma(1)}\cdots x_{(g_{\sigma(m)},h_{\sigma(m)}')\sigma(m)}. \]

Thus, \[\tilde{f}(b_1,\dots,b_m)= \sum_{\sigma\in S_m}c_\sigma \lambda_{{\bf h'}\sigma} b_{\sigma(1)}\cdots b_{\sigma(m)}.\]

On the other hand,
\begin{align*}
\phi_{\bf h}(f(w_1,\dots,w_m))
&= \phi_{\bf h}(\sum_{\sigma\in S_m}c_\sigma w_{\sigma(1)}\cdots w_{\sigma(m)})\\
&= \sum_{\sigma\in S_m}c_{\sigma}\zeta(\varphi_{\bf h}(w_{\sigma(1)}\cdots w_{\sigma(m)}))\\
&= \sum_{\sigma\in S_m}c_{\sigma}\zeta(\varphi_{\bf h}(w_{\sigma(1)})\cdots \varphi_{\bf h}(w_{\sigma(m)}))\\
&= \sum_{\sigma\in S_m}c_{\sigma}\zeta(b_{\sigma(1)}\cdots b_{\sigma(m)}) \\
&= \sum_{\sigma\in S_m}c_{\sigma}\gamma\lambda_{{\bf h'}\sigma} b_{\sigma(1)}\cdots b_{\sigma(m)}. \\
\end{align*}
The last equality follows from Lemma \ref{9}.
Now, comparing the equations the result follows. \hfill $\Box$

Finally, we state the theorem which generalizes Theorem \ref{Theorem 11}

\begin{theorem}\label{basis} 
Let $C$ be an $H$-graded color commutative superalgebra, generating the variety of all $H$-graded color commutative superalgebras and let $R$ be any $G$-graded algebra. If $\displaystyle \mathcal{E}\subseteq \bigcup_{\substack{{\bf g}\in G^n \\ n\in \mathbb{N}} } P_{\bf g}$ is a system of multilinear generators for $T_G(R)$, then the set
\[S=\{\phi_{\bf h}(f)\,|\, f\in \mathcal{E},\, {\bf h}\in H^n, n\in \mathbb{N}\}\] is a system of multilinear generators of $T_{G\times H}(R\otimes C)$.
\end{theorem}
\textit{Proof.} Let $U$ be the T-ideal generated by $S$ in $ F\langle X|G\times H \rangle$. Of course $U\subseteq T_{G\times H}(R\otimes C)$, by
Theorem \ref{3.1}. Let us now suppose that $f\in T_{G\times H}(R\otimes C)$. Since the characteristic of $K$ is zero, we may assume that $f\in P_{\bf
g\times h}$, for some ${\bf g}\in G^n$ and ${\bf h}\in H^n$. By Lemma \ref{5}, there exists $f_0\in T_G(R)\cap P_{\bf g}$ such that
$f=\phi_{\mathbf{h}}(f_0)$. Since $f_0\in T_G(R)$, $f_0\in \langle \mathcal{E}\rangle$. Then, there exist $f_1,\dots,f_n\in \mathcal{E}$,
$u_{i},v_{i}, w_j^i\in F\langle X|G\rangle$ monomials, and $\alpha_i\in F$ such that
\[f_0=\sum_i \alpha_i u_if_i(w_1^i,\dots,w_{k_i}^i)v_i\]
Since $f_0\in P_{\bf g}$, we may assume that $u_{i},v_{i}, w_j^i$ are also multilinear.  On the other hand, for each $i$, \[\phi_{\bf
h}(u_if_i(w_1^i,\dots,w_{k_i}^i)v_i)=\beta_i\phi_{\bf h}(u_i)\phi_{\bf h}(f_i(w_1^i,\dots,w_{k_i}^i))\phi_{\bf h}(v_i)\] for some $\beta_i\in F^*$.
Now Lemma \ref{10} implies that there exist ${\bf h^i}\in H^{k_i}$ and homogeneous elements $b_j^i\in F\langle X|G\times H\rangle$, such that
\[\phi_{\bf h}(f_i(w_1^i,\dots,w_{k_i}^i))=\gamma_i\phi_{\bf h^i}(f_i)(b_1^i,\dots,b_{k_i}^i),\] for some $\gamma_i\in F^*$. Hence, Theorem \ref{3.1}
implies that for each $i$, \linebreak $\phi_{\bf h^i}(f_i(b_1^i,\dots,b_{k_i}^i))\in U$, and then every summand of $\phi_{\bf h}(f_0)$ is in U. As a
consequence, $f\in U$ and $T_{G\times H}(R\otimes C)\subseteq U$. \hfill $\Box$

The above theorem has an interesting application. To show that we need to recall some results on classification of gradings on $M_n(F)$, when $F$ is an algebraic closed field of characteristic zero. Below follows a result on the classification of abelian gradings on block-triangular matrices \cite{ValentiZaicev2}, which generalize the classification of abelian gradings on full matrix algebras \cite{BahturinSegalZaicev}.

\begin{theorem}
	Let $G$ be an abelian group and let the field $F$ be algebraically closed. For any $G$-grading of the matrix algebra $UT(d_1,\dots, d_m)$ there exist integers $t$, $q_1, \dots, q_m$ such that $d_i=tq_i$, for each $i$, a subgroup $H$ of $G$ and a $q$-tuple ${\bf g}=(g_1,\dots,g_q)\in G^q$ $(q=q_1+\cdots q_m)$, such that $UT(d_1,\dots, d_m)$ is isomorphic to $M_t(F)\otimes UT(q_1,\dots, q_m)$ as a $G$-graded algebra where $M_t(F)$ is an $H$-graded algebra with a fine $H$-grading and $UT(q_1,\dots, q_m)$ has an elementary grading defined by ${\bf g}=(g_1,\dots,g_q)$.
\end{theorem}

Moreover, it turns out that such fine $H$-grading on $M_t(F)$ makes it an $H$-graded color commutative superalgebra as we can see in the next result.

\begin{theorem}[Theorem 3.4 (i), \cite{BahturinDrensky}]\label{3.4}
	Let $M_t(F)$ have a an $H$-fine grading with all homogeneous components one-dimensional. Then the $H$-graded polynomial identities $[x_{h_1},x_{h_2}]_{\beta}=0$, where $h_1$, $h_2\in H$, form a basis of the graded polynomial identities of $M_t(F)$.
\end{theorem}

The above means that $M_t(F)$ generates the variety of all $H$-graded $\beta$-commutative superalgebras.
And now we obtain

\begin{corollary}
	Let $F$ be an algebraically closed field of characteristic zero and $G$ be a finite abelian group. If the algebra $UT(d_1,\dots,d_m)$ is $G$-graded, write $UT(d_1,\dots,d_m)\cong M_t(F)\otimes UT(q_1,\dots,q_m)$ where $M_t $  has a fine $H$-grading ($H$ a subgroup of $G$), $d_i=tp_i$, for each $i$ and $UT(q_1,\dots,q_m)$ has an elementary grading induced by a $q$-tuple  $(g_1,\dots,g_q)$ of elements of $G$, with $q=q_1+\cdots+q_m$.
	
	If $g_1$, \dots, $g_q$ are pairwise distinct, then the $G\times H$-graded polynomial identities of $UT(d_1,\dots,d_m)$, follows from the polynomials $\phi_{\bf h}(f)$, with ${\bf h}\in H^n$ and $f$ in the identities of (3) -- (5) and the $G$-graded monomial identities of degree up to $2q-1$ of $UT(q_1,\dots,q_m)$.
\end{corollary}

\textit{Proof.} By Theorem  \ref{3.4}, the fine $H$-graded algebra $M_t(F)$ is an $H$-graded color commutative superalgebra generating the variety of all $H$-graded color commutative superalgebras. By Theorem \ref{Main}, a basis for the identities of $UT(q_1,\dots,q_m)$ consists of the identities (\ref{(3)}) -- (\ref{(5)}) and a finite number of graded monomial identities of degree up to $2q-1$, with $q=q_1+\cdots+q_m$. So Theorem \ref{basis} implies that the above is a basis for $UT(d_1,\dots,d_m)$, since applying the map $\phi_{\bf h}$ does not change the degree of a monomial. \hfill $\Box$

To finish this paper, we study the relation between the notion of $H$-graded color commutative superalgebras and the notion of a Regular $H$-grading.

If $H$ is a group, the notion of a regular $H$-grading was introduced by Regev and Seeman \cite{RegevSeeman} and we recall it now.

If $H$ is a group an $H$-grading on an associative algebra $A=\displaystyle\bigoplus_{h\in H} A_h$ is called a regular $H$-grading if it satisfies the following two conditions:

\noindent 1) For any $(h_1,\dots,h_n)\in H_n$, there exist $a_{h_1}\in A_{h_1}$, \dots, $a_{h_n}\in A_{h_n}$ such that $a_{h_1}\cdots a_{h_n}\neq 0$.

\noindent 2) If $g$, $h\in H$, there exist $\theta(g,h)\in F^*$ such that for any $a\in A_g$, and $b\in A_h$, $ab=\theta(g,h)ba$.

Observe that if $H$ is an abelian group, the condition 2) means that $A$ is an $H$-graded color commutative superalgebra. Indeed, we claim that map $\theta: H\times H \longrightarrow F^*$ defined above is a skew-symmetric bicharacter.

To show that, let $g,h\in H$ and consider $a\in A_g$ and $b\in A_h$ such that $ab\neq 0$.
Then on one hand one has $ab=\theta(g,h)ba$, which implies that $ba=\theta(g,h)^{-1}ab$. On the other hand $ba=\theta(h,g)ab$. Since $ab\neq 0$, we conclude that $\theta(h,g)=\theta(g,h)^{-1}$.

Also, if $g$, $h$, and $k\in H$, consider $a\in A_g$, $b\in A_h$ and $c\in A_k$ such that $abc\neq 0$. On one hand we have $abc=\theta(g,h+k)bca$. On the other hand, $abc=\theta(g,h)bac=\theta(g,h)\theta(g,k)bca$. Again, since $abc\neq 0$, and $A$ is associative, we obtain that $\theta(g,h+k)=\theta(g,h)\theta(g,k)$.

We have just proved the following result (which is also mentioned in \cite[Section 3.2]{AljadeffOfir2}):

\begin{lemma}
If an $H$-grading on an associative algebra $A$ is regular, then $A$ is a color commutative superalgebra.
\end{lemma}

Observe now that condition 1) above means that $A$ does not satisfy any monomial identity, and by condition 2) any multilinear polynomial is equivalent to a monomial. As a consequence, $A$ does not satisfies any other polynomial  identity. In particular, the polynomial identities of $A$ are all consequences of the ones obtained from condition 2), i.e., all such identities are consequence of the $\theta$-commutator. Which is the same as saying that $A$ generates the variety of $H$-graded color commutative superalgebras.

Conversely, if $A$ is an $H$-graded color commutative superalgebra generating the variety of all $H$-graded color commutative superalgebras, one easily sees that condition 1) and 2) above are verified. The above proves the next result.

\begin{theorem}
	Let $H$ be an abelian group and $A$ be an associative $H$-graded algebra. Then the $H$-grading on $A$ is regular if and only if $A$ is an $H$-graded color commutative superalgebra generating the variety of all $H$-graded color commutative superalgebras.
\end{theorem}

\begin{flushleft}
\textbf{Acknowledgements}
\end{flushleft}
This work was completed while the first author was a postdoctoral fellow at Memorial University of Newfoundland. He would like to thank prof. Yuri Bahturin for useful discussions on this subject.

\end{document}